\documentclass{amsart}

\newtheorem{thm}{Theorem}[section]

\newtheorem{remark}{Remark}[section]

\newtheorem{alphthm}{Theorem}[section]

\newtheorem{alphprop}{Proposition}[section]

\newtheorem{alphlemma}{Lemma}[section]

\newcommand{\lct}{\; \raisebox{-.96ex}{$\stackrel{\textstyle <}{\sim}$} \;}
\newcommand{\gct}{\; \raisebox{-.96ex}{$\stackrel{\textstyle >}{\sim}$} \;}
\newcommand{\mbb}{\mathbb}

\begin{document}

\title{Improved weighted restriction estimates in $\mbb R^3$}

\author{Bassam Shayya}
\address{Department of Mathematics\\
         American University of Beirut\\
         Beirut\\
         Lebanon}
\email{bshayya@aub.edu.lb}

\date{June 13, 2022}


\begin{abstract}
Suppose $0 < \alpha \leq n$, $H: \mbb R^n \to [0,1]$ is a Lebesgue 
measurable function, and $A_\alpha(H)$ is the infimum of all numbers $C$ for 
which the inequality $\int_B H(x) dx \leq C R^\alpha$ holds for all balls 
$B \subset \mbb R^n$ of radius $R \geq 1$. After Guth introduced polynomial 
partitioning to Fourier restriction theory, weighted restriction estimates 
of the form $\| Ef \|_{L^p(B,Hdx)} \lct R^\epsilon A_\alpha(H)^{1/p} 
\| f \|_{L^q(\sigma)}$ have been studied and proved in several papers, 
leading to new results about the decay properties of spherical means of 
Fourier transforms of measures and, in some cases, to progress on Falconer's 
distance set conjecture in geometric measure theory. This paper improves on 
the known estimates when $E$ is the extension operator associated with the 
unit paraboloid ${\mathcal P} \subset \mbb R^3$, reaching the full possible 
range of $p,q$ exponents (up to the sharp line) for 
$p \geq 3 + (\alpha-2)/(\alpha+1)$ and $2 < \alpha \leq 3$.
\end{abstract}

\maketitle

\section{Introduction}

This paper is concerned with the extension operator $E$ for the unit 
paraboloid
\begin{displaymath}
{\mathcal P} = \{ (\omega, |\omega|^2) : |\omega| \leq 1 \} 
\subset \mbb R^3,
\end{displaymath}
endowed with the measure $\sigma$ defined as the pushforward of Lebesgue 
measure under the map $\omega \mapsto (\omega, |\omega|^2)$ from the unit 
ball in $\mbb R^2$ to ${\mathcal P}$. The operator $E$ assigns to each 
function $f \in L^1(\sigma)$ a function $Ef: \mbb R^3 \to \mbb C$ given by
\begin{displaymath}
Ef(x) = \int e^{-2 \pi i x \cdot \xi} f(\xi) d\sigma(\xi)
      = \int_{|\omega| \leq 1} e^{-2 \pi i x \cdot (\omega, |\omega|^2)} 
	    f(\omega,|\omega|^2) d\omega.	
\end{displaymath}
One of the central problems in harmonic analysis, asks for the best possible
range of exponents $p$ and $q$ for which the estimate
\begin{equation}
\label{lebspbds}
\| Ef \|_{L^p(\mbb R^3)} \lct \| f \|_{L^q(\sigma)}
\end{equation}
holds. This problem is the three-dimensional version of the restriction 
conjecture, which (in $\mbb R^3$) asserts that (\ref{lebspbds}) is true 
whenever $p > 3$ and $4/p + 2/q \leq 2$. 

We currently know that (\ref{lebspbds}) is true whenever
\begin{equation}
\label{currentrange}
p > 3 + \frac{1}{4} \hspace{0.25in} \mbox{ and } \hspace{0.25in}
\frac{4}{p} + \frac{2}{q} \leq 2.
\end{equation}
The $p > 3 + 1/4$ and $q=\infty$ range was established by Guth in 
\cite{guth:poly} using his breakthrough polynomial partitioning method. The 
$4/p + 2/q < 2$ range was then reached in \cite{plms12046}, and the
$4/p + 2/q = 2$ range in \cite{jk:someremarks}.

When $q=\infty$, the range of $p$ in (\ref{lebspbds}) has been recently 
improved to $p > 3 + 3/13$ in \cite{42over13}.

In this paper, we use polynomial partitioning to prove weighted variants of 
(\ref{lebspbds}): we replace the $L^p(\mbb R^3)=L^p(dx)$ norm on the 
left-hand side of (\ref{lebspbds}) by the $L^p(Hdx)$ norm for an 
appropriate non-negative Lebesgue measurable function $H$ on $\mbb R^3$, and 
we obtain in return an improvement on the range of exponents in 
(\ref{currentrange}). For example, if $0 < a \leq 1$ and $H_a$ is the 
characteristic function of the set 
\begin{displaymath}
\cup_{m \in \mbb Z} \Big( [-10,10] \times \mbb R^2
+ \big( (\mbox{sgn} \, m) |m|^{1/a}, 0, 0 \big) \Big),
\end{displaymath}
where $\mbox{sgn} \, m= m/|m|$ if $m \not= 0$, and $\mbox{sgn} \, 0= 1$, 
then we prove that
\begin{equation}
\label{hsuba}
\| Ef \|_{L^p(H_a dx)} \lct \| f \|_{L^q(\sigma)}
\end{equation}
for all $f \in L^q(\sigma)$ whenever
\begin{equation}
\label{hsubarange}
p > 3 + \frac{a}{3+a} \hspace{0.25in} \mbox{ and } \hspace{0.25in}
\frac{3+a}{p} + \frac{2}{q} < 2.
\end{equation}
A simple Knapp example argument shows that the range $(3+a)/p + 2/q < 2$ in 
(\ref{currentrange}) is sharp up to the critical line $(3+a)/p + 2/q = 2$, 
in the sense that $(3+a)/p + 2/q \leq 2$ is a necessary condition for 
(\ref{hsuba}) to hold. We will show the details of this argument during the
proof of Theorem \ref{newhighdim} in Section 3. (We prove (\ref{hsuba}) in 
the second paragraph following Remark \ref{smallestp} in Section 2.)

We note that when $a=1$, (\ref{hsubarange}) becomes (\ref{currentrange}) 
(minus the sharp line $4/p + 2/q = 2$). We also note that when $q=\infty$, 
the range $p > 3 + a/(3+a)$ in (\ref{hsubarange}) improves on the range 
$p > 3 + 3/13$ of \cite{42over13} for $0 < a < 0.9$.

The best previously known range of exponents for which (\ref{hsuba}) holds 
is considerably smaller than (\ref{hsubarange}), and is the case $b=1$ (and
hence $\alpha=2+a$) of \cite[Corollary 2.1]{plms12046}. (See also Theorem
\ref{prevhighdim} and the paragraph following the statement of Theorem
\ref{newhighdim} below.)

The weighted restriction estimates that we are interested in have been 
recently studied in \cite{plms12046}, \cite{dgowwz:falconer}, and 
\cite{pems201030}. This paper uses the polynomial partitioning method from 
\cite{guth:poly} as adapted to the weighted setting in \cite{plms12046} and 
\cite{dgowwz:falconer} and obtains new and sharp (up to the critical line) 
restriction estimates (such as (\ref{hsuba})).

In addition to Guth's polynomial partitioning, \cite{dgowwz:falconer} 
employs the refined Strichartz estimates that were proved in 
\cite{dgl:schrodinger} using the Bourgain-Demeter $l^2$ decoupling theorem
\cite{bd:decoupling}. In proving our main theorem, we also follow 
\cite{dgowwz:falconer} in using the refined Strichartz estimates and combine 
this with some of the ideas of \cite{plms12046}. 

We would like to conclude this section by explaining one particular idea 
from \cite{plms12046} that we will use in this paper - continuing to use the 
weight function $H_a$ from (\ref{hsuba}) as an example - hoping to provide 
the reader with a global view of the argument.

The function $H_a$ satisfies the following dimensionality property:
\begin{equation}
\label{dimprop}
\int_{B(x_0,R)} H_a(x) dx \leq A(H_a) R^\alpha 
\end{equation}
for all $x_0 \in \mbb R^3$ and $R \geq 1$, where $\alpha=2+a$, $A(H_a)$ is
a constant that only depends on $\alpha$, and $B(x_0,R)$ is the ball of 
center $x_0$ and radius $R$. 

To prove (\ref{hsuba}) following the strategy of \cite{plms12046}, we first 
decompose the paraboloid ${\mathcal P}$ into caps $\{ \tau \}$ each of 
radius $K^{-1}$ for some large constant $K$, write $Ef=\sum_\tau Ef_\tau$
with $f_\tau = \chi_\tau f$, and use Guth's polynomial partitioning method
to bound the $L^p(B(0,R),H_a dx)$ norm of the broad part of $Ef$ 
(see Subsection 4.2) associated with the decomposition $\{ \tau \}$ of 
${\mathcal P}$. In this bound, the constant $A(H_a)$ will appear raised to 
the power $3/(\alpha+1)=3/(a+3) < 1$ (see (\ref{priortolocal}) in Subsection 
5.2).

Next, we bound the $L^p(B(0,R),H_a dx)$ norm of each $Ef_\tau$ by 
performing parabolic scaling at scale $r = K^{-1}$ and using induction on 
the radius $R$.

The parabolic scaling results in replacing the weight function $H_a$ by a 
weight $\widetilde{H}_a$ that obeys the dimensionality property 
\begin{displaymath}
\int_{B(x_0,R)} \widetilde{H}_a(x) dx \leq A(\widetilde{H}_a) R^\alpha 
\end{displaymath}
with the new constant $A(\widetilde{H}_a)$ satisfying
\begin{displaymath}
A(\widetilde{H}_a) \leq (192) r^{3-\alpha} A(H_a)
\end{displaymath}
(for a proof of this, we refer the reader to \cite[Equation (34)]{plms12046}
or \cite[Lemma 2.1]{dgowwz:falconer}). Since the bound on the broad part has 
$A(\widetilde{H}_a)$ raised to the power $3/(\alpha+1) < 1$, the power of
$r$ will consequently be reduced from $3-\alpha$ to 
$3(3-\alpha)/(\alpha+1)$.  

During the induction on the radius argument, however, one realizes that to 
reach the best possible range of $p,q$ exponents for which (\ref{hsuba}) 
holds, one needs to preserve the `original' $3-\alpha$ power of $r$.

It was realized in \cite{plms12046} that it is possible preserve the 
$3-\alpha$ power of $r$, if we take advantage of the fact that the
exponent of $R$ in (\ref{dimprop}) (which, starting from the next section, 
we will think of as the {\it dimension} of the weight) is invariant under 
localization. More precisely, if we tile $\mbb R^3$ with cubes $\{ Q_l \}$ 
each of side-length $K/3$, and define the function $H_{a,K}$ by
\begin{displaymath}
H_{a,K}(x) = A(H_a)^{-1} K^{-\alpha} \int_{Q_l} H_a(y) dy 
\hspace{0.25in} \mbox{ for } \hspace{0.25in} x \in Q_l,
\end{displaymath}
then
\begin{displaymath}
\int_{B(x_0,R)} H_{a,K}(x) dx \leq K^{3-\alpha} R^\alpha
\end{displaymath}
for all $x_0 \in \mbb R^3$ and $R \geq 1$.

So, to get the full range of exponents in (\ref{hsubarange}), we go back and 
adjust the above argument by first localizing the weight $H_a$ at scale $K$, 
as described in the previous paragraph. Then we use the locally constant 
property of the Fourier transform and the fact that $f_\tau$ is supported on 
a cap of radius $r=K^{-1}$ to see that
\begin{displaymath}
\int_{B(0,R)} |Ef_\tau(x)|^p H_a(x) dx 
\sim A(H_a) r^{3-\alpha} \int_{B(0,R)} |Ef_\tau(x)|^p H_{a,K}(x) dx. 
\end{displaymath}
Performing parabolic scaling to the integral on the right-hand side will now 
result in replacing the weight $H_{a,K}$ by a weight $\widetilde{H}_{a,K}$ 
that obeys the dimensionality property
\begin{displaymath}
\int_{B(x_0,R)} \widetilde{H}_{a,K}(x) dx 
\leq (192) r^{3-\alpha} K^{3-\alpha} R^\alpha = (192) R^\alpha
\end{displaymath}
for all $x_0 \in \mbb R^3$ and $R \geq 1$, which will allow for no losses in 
the power of $r$. (For the full argument, we refer the reader to
\cite[Section 12]{plms12046}.) 

One consequence of the above discussion is that to prove a favorable 
restriction estimate for a particular weight (such as $H_a$) that satisfies
a certain dimensionality property (such as (\ref{dimprop})), one needs to 
establish the estimate for all (properly normalized) weight functions that 
obey the same dimensionality property. This brings us to the definition with
which we start the next section.

\section{Results} 

Suppose $0 < \alpha \leq n$ and $H$ is non-negative Lebesgue measurable 
function on $\mbb R^n$ with $\| H \|_{L^\infty(\mbb R^n)} \leq 1$. Following 
\cite{plms12046} and \cite{pems201030}, we define
\begin{displaymath}
A_\alpha(H)= \inf \Big\{ C : \int_{B(x_0,R)} H(x) dx \leq C R^\alpha
\mbox{ for all } x_0 \in \mbb R^n \mbox{ and } R \geq 1 \Big\}.
\end{displaymath}
We say $H$ is a {\it weight of fractal dimension} $\alpha$ if 
$A_\alpha(H) < \infty$. (For the motivation behind referring to $\alpha$ as
a {\it fractal dimension}, see \cite[Sections 4 and 8]{pems201030}.)

Also, following \cite{dgowwz:falconer}, we denote by 
${\mathcal F}_{\alpha,n}$ the set of all non-negative Lebesgue measurable 
functions ${\mathcal H}$ on $\mbb R^n$ that satisfy
\begin{displaymath}
\int_{B(x_0,R)} {\mathcal H}(x) dx \leq R^\alpha
\end{displaymath}
for all $x_0 \in \mbb R^n$ and $R \geq 1$. 

We note that if $H$ is a weight of fractal dimension $\alpha$, then 
$A_\alpha(H)^{-1} H \in {\mathcal F}_{\alpha,n}$. Conversely, if 
${\mathcal H} \in {\mathcal F}_{\alpha,n}$ and $N \in \mbb N$, then
\begin{displaymath}
A_\alpha \big( N^{-1} \chi_{\{ {\mathcal H} \leq N \}} {\mathcal H} \big) 
\leq \frac{1}{N}.
\end{displaymath} 
Thus, if $G : \mbb R^n \to \mbb C$ is a Lebesgue measurable function, then 
the inequality
\begin{equation}
\label{unify1}
\int |G(x)| {\mathcal H}(x) dx \leq B 
\hspace{0.25in} \mbox{ for all ${\mathcal H} \in {\mathcal F}_{\alpha,n}$}
\end{equation}
is equivalent to the inequality
\begin{equation}
\label{unify2}
\int |G(x)| H(x) dx \leq A_\alpha(H) B 
\hspace{0.25in} \mbox{ for all weights $H$ of fractal dimension $\alpha$.}
\end{equation}
But, even with this equivalence, it is still more convenient to work with 
weights $H$ of fractal dimension $\alpha$ to prove the estimates that we are
interested in. In fact, for such weights, the assumption 
$\| H \|_{L^\infty(\mbb R^n)} \leq 1$ allows us to use H\"{o}lder's 
inequality to get 
\begin{eqnarray*}
\int_{B(0,R)} |G(x)| H(x) dx 
& \leq & \| G \|_{L^p(\mbb R^n)} 
         \Big( \int_{B(0,R)} H(x)^{p'} dx \Big)^{1/p'} \\
& \leq & \| G \|_{L^p(\mbb R^n)} \Big( \int_{B(0,R)} H(x) dx \Big)^{1/p'} \\		 
& \leq & A_\alpha(H)^{1/p'} R^{\alpha/p'} \| G \|_{L^p(\mbb R^n)} 
\end{eqnarray*}
for all $1 < p < \infty$. For ${\mathcal H} \in {\mathcal F}_{\alpha,n}$, 
however, the absence of a bound on $\| {\mathcal H} \|_{L^\infty(\mbb R^n)}$ 
prevents us from carrying out such calculations so easily (except if $G$ has 
compact Fourier support; see \cite[Remark 1.3]{dgowwz:falconer}).

In proving the main result of this paper (Theorem \ref{mainjj} below), the 
losses incurred as a result of applying H\"{o}lder's inequality as above 
will be compensated for by the localization argument that was discussed in 
Section 1 (see also Subsection 5.3).

The equivalence between {\rm (\ref{unify1})} and {\rm (\ref{unify2})} will
allow us to unify the notation of \cite{plms12046}, \cite{dgowwz:falconer}, 
and \cite{pems201030} and present their weighted restriction estimates (in 
$\mbb R^3$) in the following three theorems; the first of which describes 
the current state of affairs in low fractal dimensions. 

\begin{alphthm}[{\cite[Theorem 2.1]{pems201030}}]
\label{prevlowdim}
{\rm (i)} Suppose $0 < \alpha \leq 1$ and $p > 2$. Then
\begin{displaymath}
\| Ef \|_{L^p(H dx)} \lct A_\alpha(H)^{1/p} \| f \|_{L^2(\sigma)}
\end{displaymath}
for all $f \in L^2(\sigma)$ and weights $H$ of fractal dimension $\alpha$.

{\rm (ii)} Suppose $1 < \alpha \leq 3/2$ and $p > 2 \alpha$. Then
\begin{displaymath}
\| Ef \|_{L^p(H dx)} \lct A_\alpha(H)^{1/p} \| f \|_{L^2(\sigma)}
\end{displaymath}
for all $f \in L^2(\sigma)$ and weights $H$ of fractal dimension $\alpha$.
\end{alphthm}

In intermediate fractal dimensions in $\mbb R^3$, the estimates we know are 
local in the sense that the $L^p(H dx)$ norm of $Ef$ is taken over a ball of
radius $R$ rather than over the entire $\mbb R^3$. For some weights (such as 
$H_a$ in (\ref{hsuba})) the local estimates can be turned into global ones 
by using Tao's $\epsilon$-removal lemma from \cite{tt:removal}. For the 
details of this argument, we refer the reader to 
\cite[Proof of Corollary 2.1]{plms12046}.

\begin{alphthm}[{\cite[Theorem 1.4]{dgowwz:falconer}}]
\label{prevousdim}
Suppose $3/2 < \alpha \leq 2$. Then to every $\epsilon > 0$ there is a 
constant $C_\epsilon$ such that
\begin{displaymath}
\| Ef \|_{L^3(B(0,R), H dx)} 
\leq C_\epsilon R^\epsilon A_\alpha(H)^{1/3} \| f \|_{L^2(\sigma)}  
\end{displaymath}
for all $f \in L^2(\sigma)$, weights $H$ of fractal dimension $\alpha$, and 
$R \geq 1$.
\end{alphthm}

In high fractal dimensions, i.e.\ $2 < \alpha \leq 3$, the best previously 
known results are stated in the following theorem. 

\begin{alphthm}[{\cite[Theorem 5.1]{plms12046}}]
\label{prevhighdim}
{\rm (i)} Suppose $2 < \alpha < 5/2$ and
\begin{displaymath}
p=3 + \frac{2\alpha-3}{2\alpha+3} 
\hspace{0.25in} \mbox{ and } \hspace{0.25in} 
q = \frac{2p}{3}. 
\end{displaymath}
Then to every $\epsilon > 0$ there is a constant $C_\epsilon$ such that
\begin{displaymath}
\| Ef \|_{L^p(B(0,R), H dx)} 
\leq C_\epsilon R^\epsilon A_\alpha(H)^{1/p} \| f \|_{L^q(\sigma)} 
\end{displaymath}
for all $f \in L^q(\sigma)$, weights $H$ of fractal dimension $\alpha$, and 
$R \geq 1$.

{\rm (ii)} Suppose $5/2 \leq \alpha \leq 3$, $p=13/4$, and 
$q > (11-2\alpha)/13$. Then to every $\epsilon > 0$ there is a constant 
$C_\epsilon$ such that
\begin{displaymath}
\| Ef \|_{L^p(B(0,R), H dx)} 
\leq C_\epsilon R^\epsilon A_\alpha(H)^{1/p} \| f \|_{L^q(\sigma)} 
\end{displaymath}
for all $f \in L^q(\sigma)$, weights $H$ of fractal dimension $\alpha$, and 
$R \geq 1$.
\end{alphthm}

We alert the reader that in \cite[Theorem 5.1]{plms12046}, the restriction
estimates are written in dual form: from $L^{p'}(B(0,R), H dx)$ to 
$L^{q'}(\sigma)$ (see also Section 3 below). Also, part (i) of Theorem 
\ref{prevhighdim} was in fact proved in \cite{plms12046} for 
$3/2 \leq \alpha < 5/2$, but for $\alpha=3/2$ its result has been surpassed 
by \cite[Theorem 2.1]{pems201030} and for $3/2 < \alpha \leq 2$ by 
\cite[Theorem 1.4]{dgowwz:falconer} as presented in Theorems 
\ref{prevlowdim} and \ref{prevousdim} above.

The aim of this paper is to prove the following theorem.

\begin{thm}
\label{mainjj}
Suppose $2 < \alpha \leq 3$, $p=(4\alpha+1)/(\alpha+1)$, and 
$2 \leq \gamma < 2p-\alpha-1$. Then to every $\epsilon > 0$ there is a 
constant $C_\epsilon$ such that
\begin{displaymath}
\int_{B(0,R)} |Ef(x)|^p H(x) dx \leq C_\epsilon R^\epsilon A_\alpha(H)
\| f \|_{L^2(\sigma)}^\gamma \| f \|_{L^\infty(\sigma)}^{p-\gamma} 
\end{displaymath}
for all $f \in L^\infty(\sigma)$, weights $H$ of dimension $\alpha$, and 
$R \geq 1$.
\end{thm}

Theorem \ref{mainjj} will allow us to extend the range of exponents for 
which our weighted estimates hold in high fractal dimensions until 
(but still not including) the sharp line $(\alpha+1)/p + 2/q = 2$, and hence 
replace Theorem \ref{prevhighdim} by the following result.

\begin{thm}
\label{newhighdim}
Suppose $2 < \alpha \leq 3$. Then:
\\
{\rm (i)} Suppose
\begin{displaymath}
p \geq 3 + \frac{\alpha-2}{\alpha+1} \hspace{0.25in} \mbox{ and } 
\hspace{0.25in} \frac{\alpha+1}{p} + \frac{2}{q} < 2.
\end{displaymath} 
Then to every $\epsilon > 0$ there is a constant $C_\epsilon$ such that
\begin{displaymath}
\| Ef \|_{L^p(B(0,R), H dx)} 
\leq C_\epsilon R^\epsilon A_\alpha(H)^{1/p} \| f \|_{L^q(\sigma)} 
\end{displaymath}
for all $f \in L^q(\sigma)$, weights $H$ of fractal dimension $\alpha$, and 
$R \geq 1$.
\\
{\rm (ii)} Suppose $1 \leq p,q \leq \infty$ are a pair of exponents for 
which the estimate in part {\rm (i)} holds. Then
\begin{displaymath}
p \geq \alpha \hspace{0.25in} \mbox{ and } \hspace{0.25in}
\frac{\alpha+1}{p} + \frac{2}{q} \leq 2.
\end{displaymath}
\end{thm}

\begin{remark}
\label{smallestp}
The lower bound $p \geq \alpha$ in part {\rm (ii)} of Theorem 
\ref{newhighdim} was proved by the author in \cite[Theorem 2.3]{pems201030}. 
We are restating it here for obvious reasons.
\end{remark}

We will prove Theorem \ref{newhighdim} (using Theorem \ref{mainjj}) in the 
next section. The proof of Theorem \ref{mainjj} will occupy all the sections
of the paper following the next one.

Going back to (\ref{hsuba}), we know from (\ref{dimprop}) that $H_a$ is a
weight of fractal dimension $\alpha=2+a$ with $A_\alpha(H) \lct 1$, so 
applying part (i) of Theorem \ref{newhighdim} and using Tao's 
$\epsilon$-removal lemma (as discussed in the paragraph preceding the 
statement of Theorem \ref{prevousdim} above), we see that (\ref{hsuba}) 
holds for all pairs $p,q$ of exponents that satisfy (\ref{hsubarange}).

We would like to point out to the reader that 
\cite[Theorem 1.4]{dgowwz:falconer} also gives an estimate in the regime 
$2 < \alpha \leq 3$: 
\begin{displaymath}
\| Ef \|_{L^3(B(0,R), {\mathcal H} dx)} 
\leq C_\epsilon R^\epsilon R^{(\alpha-2)/3} \| f \|_{L^2(\sigma)} 
\end{displaymath}
for all $f \in L^2(\sigma)$, ${\mathcal H} \in {\mathcal F}_{\alpha,3}$, and 
$R \geq 1$, which the authors of \cite{dgowwz:falconer} then use to obtain 
decay estimates on the $L^2$ spherical means of the Fourier transforms of 
measures (see \cite[Theorem 1.6]{dgowwz:falconer}) that were new at the time 
of writing of their paper. Theorem \ref{newhighdim} improves on the decay 
estimates of \cite{dgowwz:falconer} in the regime 
$2 < \alpha < 1+(2/\sqrt{3})$. But all these decay estimates are now 
inferior to those of \cite{dz:schrodinger3}, so we omit the details.

\section{Proof of Theorem \ref{newhighdim}}

(i) We will establish the estimate in dual form. For $g \in C_0(B(0,R))$, we
let ${\mathcal R} g = \widehat{gH}$, and our goal is to show that
\begin{equation}
\label{dualform}
\| {\mathcal R} g \|_{L^{q'}(\sigma)} 
\lct R^\epsilon A_\alpha(H)^{1/p} \| g \|_{L^{p'}(Hdx)}.
\end{equation}
To do this, we will use Theorem \ref{mainjj} in the following way:
\begin{eqnarray*}
\Big| \int {\mathcal R} g(\xi) f(\xi) d\sigma(\xi) \Big| 
&  =   & \Big| \int Ef(x) g(x) H(x) dx \Big| \\
& \leq & \| Ef \|_{L^p(B(0,R), Hdx)} \| g \|_{L^{p'}(Hdx)} \\
& \lct & R^\epsilon A_\alpha(H)^{1/p} \| f \|_{L^2(\sigma)}^{\gamma/p}	 
	     \| f \|_{L^\infty(\sigma)}^{(p-\gamma)/p} \| g \|_{L^{p'}(Hdx)}
\end{eqnarray*}
for all $f \in L^1(\sigma)$. The key idea now (which the author has 
previously used in \cite[Theorem 5.1]{plms12046}) is that the above estimate
allows us to estimate the $\sigma$-measure of the superlevel sets 
$\{ \xi \in {\mathcal P} : |{\mathcal R} g(\xi)| > \lambda \}$ for 
$0 < \lambda \leq \| g \|_{L^1(Hdx)}$ as follows. 

Letting
\begin{displaymath}
S_{\lambda,l} = \{ \xi \in {\mathcal P} : 
2^{l-1} \lambda < |{\mathcal R}g(\xi)| \leq 2^l \lambda \}
\end{displaymath}
and inserting $\overline{{\mathcal R}g} \; \chi_{S_{\lambda,l}}$ for $f$ in 
the above estimate, we see that
\begin{displaymath}
\int_{S_{\lambda,l}} |{\mathcal R}g(\xi)|^2 d\sigma(\xi) 
\lct R^\epsilon A_\alpha(H)^{1/p} 
\| {\mathcal R}g \|_{L^2(\chi_{S_{\lambda,l}}d\sigma)}^{\gamma/p}
(2^l \lambda)^{(p-\gamma)/p} \| g \|_{L^{p'}(Hdx)},
\end{displaymath}
so that
\begin{displaymath}
(2^l \lambda)^2 \sigma(S_{\lambda,l}) 
\lct R^\epsilon A_\alpha(H)^{1/p} (2^l \lambda)^{\gamma/p} 
     \sigma(S_{\lambda,l})^{\gamma/(2p)} 
	 (2^l \lambda)^{(p-\gamma)/p} \| g \|_{L^{p'}(Hdx)},
\end{displaymath}
so that
\begin{displaymath}
(2^l \lambda)^{2p} \sigma(S_{\lambda,l})^{2p-\gamma} 
\lct R^{2p\epsilon} A_\alpha(H)^2 \| g \|_{L^{p'}(Hdx)}^{2p}.
\end{displaymath}
Letting $q'=2p/(2p-\gamma)$, this becomes
\begin{displaymath}
\sigma(S_{\lambda,l}) \lct R^{q'\epsilon} A_\alpha(H)^{q'/p}
\| g \|_{L^{p'}(Hdx)}^{q'} (2^l \lambda)^{-q'}.
\end{displaymath}
Thus
\begin{displaymath}
\sigma(\{ \xi \in {\mathcal P} : |{\mathcal R}g(\xi)| > \lambda \})
\leq \sum_{l=1}^\infty \sigma(S_{\lambda,l})
\lct R^{q'\epsilon} A_\alpha(H)^{q'/p} \| g \|_{L^{p'}(Hdx)}^{q'} 
     \lambda^{-q'}.
\end{displaymath}
Of course, we also have the trivial bound $\sigma(\{ \xi \in {\mathcal P} : 
|{\mathcal R}g(\xi)| > \lambda \}) \leq \sigma({\mathcal P})$, so
\begin{equation}
\label{suplevelbd}
\sigma(\{ \xi \in {\mathcal P} : |{\mathcal R}g(\xi)| > \lambda \}) \leq 
\min \Big[ \sigma({\mathcal P}), R^{q'\epsilon} A_\alpha(H)^{q'/p} 
           \| g \|_{L^{p'}(Hdx)}^{q'} \lambda^{-q'} \Big].
\end{equation}

Having obtained a good bound on the size of the superlevel sets of 
${\mathcal R}g$, we can now estimate its $L^{q'}(\sigma)$ norm:
\begin{displaymath}
\| {\mathcal R}g \|_{L^{q'}(\sigma)}^{q'} = \int_0^{\lambda_1} 
\sigma(\{ \xi \in {\mathcal P} : |{\mathcal R}g(\xi) > \lambda \}) \,
\lambda^{q'-1} d\lambda,
\end{displaymath} 
where $\lambda_1 = \| g \|_{L^1(Hdx)}$. We also let
\begin{displaymath}
\lambda_0= R^\epsilon A_\alpha(H)^{1/p} \| g \|_{L^{p'}(Hdx)}
\end{displaymath} 

If $\lambda_0 \geq \lambda_1$, then the trivial part of (\ref{suplevelbd})
tells us that
\begin{displaymath}
\| {\mathcal R}g \|_{L^{q'}(\sigma)}^{q'} \leq \int_0^{\lambda_0} 
\sigma({\mathcal P}) \, \lambda^{q'-1} d\lambda \lct \lambda_0^{q'},
\end{displaymath}
so that
\begin{displaymath}
\| {\mathcal R}g \|_{L^{q'}(\sigma)} 
\lct R^\epsilon A_\alpha(H)^{1/p} \| g \|_{L^{p'}(Hdx)}.
\end{displaymath}

If $\lambda_0 \leq \lambda_1$, then, in view of the inequality before the 
last, we only have to estimate
\begin{displaymath}
\int_{\lambda_0}^{\lambda_1} 
\sigma(\{ \xi \in {\mathcal P} : |{\mathcal R}g(\xi) > \lambda \}) \, 
\lambda^{q'-1} d\lambda,
\end{displaymath}
which, by the non-trivial part of (\ref{suplevelbd}), is
\begin{eqnarray*}
\lct \int_{\lambda_0}^{\lambda_1}
        R^{q'\epsilon} A_\alpha(H)^{q'/p} \| g \|_{L^{p'}(Hdx)}^{q'} 
		\lambda^{-1} d\lambda
= \Big( \log \frac{\lambda_1}{\lambda_0} \Big) 
  R^{q'\epsilon} A_\alpha(H)^{q'/p} \| g \|_{L^{p'}(Hdx)}^{q'}.
\end{eqnarray*}
By H\"{o}lder's inequality (applied with the measure $Hdx$ on the ball 
$B(0,R)$),
\begin{displaymath}
\lambda_1 \leq \Big( \int_{B(0,R)} H(x) dx \Big)^{1/p} \| g \|_{L^{p'}(Hdx)}  
\leq A_\alpha(H)^{1/p} R^{\alpha/p} \| g \|_{L^{p'}(Hdx)},
\end{displaymath}
so
\begin{displaymath}
\log \frac{\lambda_1}{\lambda_0}  
\leq \log \frac{A_\alpha(H)^{1/p} R^{\alpha/p} \| g \|_{L^{p'}(Hdx)}}
               {R^\epsilon A_\alpha(H)^{1/p} \| g \|_{L^{p'}(Hdx)}}
\lct \log R,
\end{displaymath}
and so
\begin{displaymath}
\| {\mathcal R}g \|_{L^{q'}(\sigma)} 
\lct R^\epsilon A_\alpha(H)^{1/p} \| g \|_{L^{p'}(Hdx)}.
\end{displaymath}

Recalling that $q'=2p/(2p-\gamma)$ and $2 \leq \gamma <2p-\alpha-1$, we see
that (\ref{dualform}) is true whenever $1/p \leq 1/q < 1-(\alpha+1)/(2p)$, 
and hence (by H\"{o}lder) whenever $1/q < 1-(\alpha+1)/(2p)$, as promised.

(ii) As Remark \ref{smallestp} tells us, the lower bound $p \geq \alpha$ on 
the exponent $p$ comes from \cite[Theorem 2.3]{pems201030}, so we only need 
to show that $(\alpha+1)/p + 2/q \leq 2$. For this we use a standard
Knapp-example argument.

To every $R > 1$ there is a function $f_R$ on ${\mathcal P}$ such that 
$|f_R|=1$ on the cap $\{ (\omega,|\omega|^2) : |\omega| \leq R^{-1/2} \}$ 
and $|Ef_R| \gct R^{-1}$ on the box 
$[-R^{1/2},R^{1/2}] \times [-R^{1/2},R^{1/2}] \times [-R,R]$.

If the estimate in part (i) of the theorem holds for all weights $H$ of 
fractal dimension $\alpha$, then, in particular, it will hold for the weight 
$H_a$ from Section 1 with $a = \alpha-2 \in (0,1]$. Since 
$|m|^{1/a} \leq R^{1/2}$ for $\sim R^{a/2}$ integers $m$, we have
\begin{displaymath}
\int_{B(0,10R)} |Ef_R(x)|^p H_a(x) dx \gct R^{-p} R^{a/2} R^{3/2} 
= R^{-p} R^{(\alpha+1)/2}.
\end{displaymath} 
Also, $\| f_R \|_{L^q(\sigma)} \sim R^{-1/q}$, so
$R^{-p} R^{(\alpha+1)/2} \lct R^\epsilon R^{-p/q}$ for all sufficiently 
large $R$, and so $(\alpha+1)/p + 2/q \leq 2$.

\section{Overview of Guth's polynomial partitioning method}

If $B$ is a ball in $\mbb R^2$, then the {\it cap corresponding to} $B$ is 
defined to be the set $\Theta \subset \mbb R^3$ given by 
$\Theta = \{ (\omega,|\omega|^2) : \omega \in B^2 \}$. If $B$ has center 
$\omega_0$ and radius $\rho$, then $\Theta$ has center 
$\xi_0= (\omega_0,|\omega_0|^2)$ and radius $\rho$. Moreover, if $C$ is a 
constant, then $C \Theta$ is the cap of the same center as $\Theta$ and $C$ 
times its radius.

Let $m$ and $r$ be positive numbers and $\{ \Theta \}$ be caps corresponding 
to a collection of balls that cover the unit ball 
$\mbb B^2 \subset \mbb R^2$. If the centers of the balls are in $\mbb B^2$ 
and are $r$-separated, and the radius of each ball lies in the interval 
$[r,r\sqrt{m}]$, we say that $\{ \Theta \}$ is a decomposition of the unit
paraboloid ${\mathcal P}$ into caps of radius $\sim r$ and multiplicity at 
most $m$. (It is easy to see that if a point belongs to $M$ of the caps 
$\Theta$, then $M \leq c \, m$ for some absolute constant $c$.)

One reason polynomial partitioning works well for the restriction problem is 
the wave packet decomposition of $Ef$ and the way the wave packets interact 
with zero sets of polynomials.

\subsection{The wave packet decomposition}

The wave packet decomposition was first used to prove restriction estimates
in Bourgain's paper \cite{jb:besicovitch}, and has since become a standard 
tool in Fourier restriction theory. The version that we are going to use in 
this paper comes from \cite[Section 2]{guth:poly} (see also 
\cite[Section 3]{guth:poly2} and \cite[Section 7]{plms12046}).

\begin{alphprop}
\label{wave}
Suppose $\delta > 0$, $R > 1$, $N$ is a positive integer, $\theta$ is a cap 
in ${\mathcal P}$ of center $\xi_0$ and radius $r=R^{-1/2}$, and $v(\theta)$ 
is the unit normal vector of ${\mathcal P}$ at $\xi_0$.

Then there is a countable collection $\widetilde{\mbb T}(\theta)=\{ T \}$ of
finitely overlapping tubes in $\mbb R^3$ of radius $R^{(1/2)+\delta}$, which 
are parallel to $v(\theta)$, such that the following holds. To every 
function $f \in L^2(\sigma)$ with $\mbox{\rm supp}\, f \subset \theta$ there 
is a sequence $\{ f_T \}$ in $L^2(\sigma)$ with the following properties.
\\
{\rm (i)} Each $f_T$ is supported in $3\theta$,
$f= \sum_{T \in \widetilde{\mbb T}(\theta)} f_T$ in $L^1(\sigma)$, and
\begin{displaymath}
\sum_{T \in \widetilde{\mbb T}(\theta)} \int |f_T|^2 d\sigma
\lct \| f \|_{L^2(\sigma)}^2.
\end{displaymath}
{\rm (ii)} We have
\begin{displaymath}
\sum_{T \in \widetilde{\mbb T}(\theta) : \, x \not\in T} |E f_T(x)|
\lct R^{-N} \| f \|_{L^1(S)}
\end{displaymath}
for all $x \in \mbb R^3$ with $|x| \leq R$, where the implicit constant 
depends only on $N$.
\\
{\rm (iii)} If $T_1, T_2 \in \widetilde{\mbb T}(\theta)$ are disjoint, then
\begin{displaymath}
\Big| \int f_{T_1} \overline{f_{T_2}} \; d\sigma \Big|
\lct R^{-N} \| f \|_{L^1(\sigma)}^2,
\end{displaymath}
where the implicit constant depends only on $N$.
\\
{\rm (iv)} Let $\mbb T(\theta)=\{ T \in \widetilde{\mbb T}(\theta) : T \cap
B(0,R) \not=\emptyset \}$. Then
\begin{displaymath}
\Big| Ef(x) - \sum_{T \in \mbb T(\theta)} Ef_T(x) \Big|
\lct R^{-N} \| f \|_{L^1(\sigma)}
\end{displaymath}
for all $x \in B(0,R)$.
\end{alphprop}

We note that part (iv) of Proposition \ref{wave} is an immediate consequence
of parts (i) and (ii).

We fix a decomposition $\{ \theta \}$ of the paraboloid ${\mathcal P}$ into 
caps of radius $\sim R^{-1/2}$ and multiplicity at most $1$, as defined at 
the start of Section 4. 

Any function $f$ on ${\mathcal P}$ can now be written as a sum 
$f = \sum_\theta f_\theta$ with $f_\theta$ supported in $\theta$ and such 
that the supports of $f_\theta$ and $f_{\theta'}$ are disjoint whenever
$\theta \not= \theta'$. Also, if $f \in L^2(\sigma)$, then so are the 
functions $f_\theta$. Applying Proposition \ref{wave} to each $f_\theta$, we
get that
\begin{displaymath}
f = \sum_\theta \sum_{T \in \widetilde{\mbb T}(\theta)} (f_{\theta})_T
\end{displaymath}
with the sum converging in $L^1(\sigma)$. (Since the decomposition 
$\{ \theta \}$ of ${\mathcal P}$ will be fixed for the rest of the paper, we 
will write $f_T$ for $(f_{\theta})_T$ to simplify the notation.)

Letting $\widetilde{\mbb T} = \cup_\theta \widetilde{\mbb T}(\theta)$, we 
then see that
\begin{equation}
\label{infbutuni}
Ef= \sum_{T \in \widetilde{\mbb T}} Ef_T
\end{equation}
with the sum converging uniformly. 

Letting $\mbb T = \{ T \in \widetilde{\mbb T}: T \cap B(0,R) \not= \emptyset 
\} = \cup_\theta \mbb T(\theta)$, part (iv) of Proposition \ref{wave} tells 
us that, on $B(0,R)$, we can replace the infinite sum in (\ref{infbutuni}) 
by a finite sum: 
\begin{equation}
\label{finbuterr}
Ef(x) = 
\sum_{T \in \mbb T} Ef_T(x) + O \big( R^{-N} \| f \|_{L^1(\sigma)} \big)
\end{equation}
for all $x \in B(0,R)$. This representation of $Ef$ is often referred to as 
the wave packet decomposition of $Ef$ on the ball $B(0,R)$. Part (ii) of 
Proposition \ref{wave} tells us that each wave packet $Ef_T $ (which is 
equal to $E(f_{\theta})_T$ for some cap $\theta$ that comes from our fixed 
decomposition of the unit paraboloid) is essentially supported in 
$B(0,R) \cap T$, which is a tube in $\mbb R^3$ of radius $R^{(1/2)+\delta}$ 
and length $\sim R$ and points in the direction $v(\theta)$ (since 
$T \in \mbb T(\theta)$).

Suppose $\{ \theta_1, \ldots, \theta_L \}$ is a subset of our set 
$\{ \theta \}$ of caps, and, for each $l$, ${\mathcal T}_l$ is a subset of 
the set $\mbb T(\theta_l)$ of tubes corresponding to the cap $\theta_l$. 
For $1 \leq l \leq L$, part (i) of Proposition \ref{wave} tells us that
\begin{displaymath}
\sum_{T \in {\mathcal T}_l} \int |f_T|^2 d\sigma 
\leq \sum_{T \in \widetilde{\mbb T}(\theta_l)} \int |f_T|^2 d\sigma 
\lct \int_{\theta_l} |f|^2 d\sigma
\end{displaymath}
(recall that for $T \in \widetilde{\mbb T}(\theta_l)$, 
$f_T=(f_{\theta_l})_T$). Also, writing
\begin{displaymath}
\Big| \sum_{T \in {\mathcal T}_l} f_T \Big|^2 =
\sum_{T_1 \cap T_2 = \emptyset} f_{T_1} \overline{f_{T_2}} +
\sum_{T_1 \cap T_2 \not= \emptyset} f_{T_1} \overline{f_{T_2}},
\end{displaymath}
where $(T_1,T_2)$ runs over ${\mathcal T} \times {\mathcal T}$, and using 
part (iii) of Proposition \ref{wave} to get
\begin{displaymath}
\sum_{T_1 \cap T_2 = \emptyset} 
\Big| \int f_{T_1} \overline{f_{T_2}} d\sigma \Big|
\lct \int_{\theta_l} |f|^2 d\sigma,
\end{displaymath}
and the fact that 
\begin{displaymath}
|\int f_{T_1} \overline{f_{T_2}} d\sigma| 
\leq \frac{1}{2} \int |f_{T_1}|^2 d\sigma 
     + \frac{1}{2} \int |f_{T_2}|^2 d\sigma
\end{displaymath}
to get
\begin{displaymath}
\sum_{T_1 \cap T_2 \not= \emptyset} 
\Big| \int f_{T_1} \overline{f_{T_2}} d\sigma \Big|
\lct \sum_{T \in {\mathcal T}_l} \int |f_T|^2 d\sigma,
\end{displaymath}
we see that
\begin{displaymath}
\int \Big| \sum_{T \in {\mathcal T}_l} f_T \Big|^2 d\sigma 
\lct \int_{\theta_l} |f|^2 d\sigma.
\end{displaymath}
Since for each $\theta$, the set 
$J_\theta=\{ 1 \leq l \leq L: (3\theta) \cap (3\theta_l) \not= \emptyset \}$ 
has cardinality $\!\!\! \lct 1$, it follows that
\begin{displaymath}
\int_{3\theta} \Big| \sum_{l=1}^L \sum_{T \in {\mathcal T}_l} f_T \Big|^2 
d\sigma \lct \sum_{l \in J_\theta} 
\int \Big| \sum_{T \in {\mathcal T}_l} f_T \Big|^2 d\sigma
\lct \sum_{l \in J_\theta} \int_{\theta_l} |f|^2 d\sigma 
\end{displaymath}
(recall from part (i) of Proposition \ref{wave} that if 
$T \in \mbb T(\theta_l) \subset \widetilde{\mbb T}(\theta_l)$, then $f_T$ is 
supported in $3\theta_l$). But $l \in J_\theta$ implies 
$\theta_l \subset 10 \, \theta$, so
\begin{equation}
\label{fewcaps}
\int_{3\theta} \Big| \sum_{l=1}^L \sum_{T \in {\mathcal T}_l} f_T \Big|^2 
d\sigma \lct \int_{10\theta} |f|^2 d\sigma.
\end{equation}

\subsection{The broad part of $Ef$}

The main idea introduced by Bourgain and Guth in \cite{jb:besicovitch} 
concerning the restriction problem in $\mbb R^3$ was to break down $Ef$ into 
a broad part and a narrow part, estimate the broad by combining the 
Bennett-Carbery-Tao multilinear restriction theorem from 
\cite{bct:multilinear} with Wolff's Kakeya result from \cite{tw:hairbrush} 
and the two-dimensional bilinear restriction theorem (see 
\cite{taovv:bilinear}), and estimate the narrow part by parabolic rescaling 
and induction. In \cite{guth:poly}, Guth replaced the multilinear theorem by
polynomial partitioning and Wolff's Kakeya result by a bound on the number
of the caps $\{ \theta \}$ such that $\mbb T(\theta)$ has at least one tube
that intersects the zero set of the polynomial tangentially (see part (iii) 
of Proposition \ref{tangtrans} below) that he proved by adapting Wolff's 
hairbrush argument from \cite{tw:hairbrush} to the polynomial partitioning
setting. In this subsection, we introduce the broad part of $Ef$ as defined
in \cite{guth:poly} and show how it connects to bilinear expressions.

Let $K$ be a large constant and $\{ \tau \}$ a decomposition of the unit 
paraboloid ${\mathcal P}$ into caps of radius $\sim K^{-1}$ and multiplicity 
at most $m$. Any function $f : {\mathcal P} \to \mbb C$ can then be written 
as $f=\sum_\tau f_\tau$ with $\mbox{supp} \, f_\tau \subset \tau$ (but we do
not insist this time that
$(\mbox{supp} \, f_\tau) \cap (\mbox{supp} \, f_{\tau'})= \emptyset$ if
$\tau \not=\tau'$ as we did with the decomposition $f= \sum_\theta f_\theta$
of the previous subsection).

Now let $0 < \beta \leq 1$ and $f \in L^1(\sigma)$. Given a point 
$x \in \mbb R^3$, we say that $x$ is $\beta$-broad for $Ef$ if no single 
$Ef_\tau$ dominates the value of $Ef$ at $x$. More precisely, $x$ is {\it 
$\beta$-broad} for $Ef$ if
\begin{displaymath}
\max_\tau |Ef_\tau(x)| < \beta |Ef(x)|.
\end{displaymath} 
The {\it broad part of $Ef$} is the function $\mbox{Br}_\beta Ef$ defined as
\begin{displaymath}
\mbox{Br}_\beta Ef(x) = \left\{ 
\begin{array}{ll}
|Ef(x)| & \mbox{ if $x$ is $\beta$-broad for $Ef$,} \\
0       & \mbox{ otherwise.}
\end{array} \right.
\end{displaymath}
We note that
\begin{displaymath}
|Ef(x)|^p \leq \mbox{Br}_\beta Ef(x)^p + \beta^{-p} \sum_\tau |Ef_\tau(x)|^p
\end{displaymath}
for all $x \in \mbb R^3$, so that
\begin{eqnarray*}
\lefteqn{\int_{B(0,R)} |Ef(x)|^p H(x) dx} \\
& \leq & \int_{B(0,R)} \mbox{Br}_\beta Ef(x)^p H(x) dx + \beta^{-p}
         \sum_\tau \int_{B(0,R)} |Ef_\tau(x)|^p H(x) dx.
\end{eqnarray*}
In broad terms, Guth's strategy consists of using
\begin{itemize} 
\item polynomial partitioning to upgrade the two-dimensional bilinear 
      restriction theorem into an estimate on 
      $\int_{B(0,R)} \mbox{Br}_\beta Ef(x)^p H(x) dx$
\item parabolic rescaling and induction on the radius $R$ to estimate \\ 
      $\int_{B(0,R)} |Ef_\tau(x)|^p H(x)dx$.
\end{itemize}

It will be convenient at this time to introduce the following notation. For
$(R,K,m) \in [1,\infty) \times [1,\infty) \times [1,\infty)$, we let 
$\Lambda(R,K,m)$ be the set of all functions $f \in L^1(\sigma)$ such that 
$f= \sum_\tau f_\tau$ for some decomposition $\{ \tau \}$ of ${\mathcal P}$ 
of multiplicity $m$ with $\mbox{supp} \, f_\tau \subset \tau$ and
\begin{equation}
\label{baverages}
\int_{B(\xi_0,R^{-1/2}) \cap {\mathcal P}} |f_\tau(\xi)|^2 d\sigma(\xi)
\leq \frac{1}{R}
\end{equation}
for all $\xi_0 \in {\mathcal P}$. Since ${\mathcal P}$ can be covered by 
$\sim R$ of the balls $B(\xi_0,R^{-1/2})$ that appear in (\ref{baverages}), 
it follows that
\begin{displaymath}
\int |f_\tau(\xi)|^2 d\sigma(\xi) \lct 1
\end{displaymath}
for each cap $\tau$.

(A first glance at Theorem \ref{mainjj} suggests that the natural way to 
normalize its estimate is to assume that 
$\| f \|_{L^\infty(\sigma)} \leq 1$. This condition is replaced by 
(\ref{baverages}) because we are unable to control the $L^\infty$ norms 
of the functions $f_{\tau,T}$ that arise in the wave packet decomposition of 
$Ef_{\tau,T}$.)

\subsection{Polynomial partitioning}

We denote by $\mbox{Poly}_D(\mbb R^n)$ the space of polynomials in $n$ real
variables of degree at most $D$. If $P \in \mbox{Poly}_D(\mbb R^n)$, then we
denote by $Z(P)$ the zero set of $P$. We say that $P$ is non-singular if
$\nabla P(x) \not= 0$ for all $x \in Z(P)$. 

Our starting point is the following result of Guth.

\begin{alphthm}[{\cite[Corollary 1.7]{guth:poly}}]
\label{polypart}
To every non-negative function $F \in L^1(\mbb R^n)$, with 
$\| F \|_{L^1(\mbb R^n)} > 0$, and integer $D \in \mbb N$ there is a 
polynomial $P \in \mbox{\rm Poly}_D(\mbb R^n) \setminus \{ 0 \}$, which is a
product of non-singular polynomials, and a positive constant $C_n$, which
only depends on the dimension $n$, such that:
\\
{\rm (i)} $\mbb R^n \setminus Z(P)$ is a disjoint union of open sets $O_i$.
\\
{\rm (ii)} For each $i$, we have
\begin{displaymath}
C_n^{-1} D^{-n} \| F \|_{L^1(\mbb R^n)} \leq \int_{O_i} F(x) dx 
\leq C_n D^{-n} \| F \|_{L^1(\mbb R^n)}.
\end{displaymath}
\end{alphthm}

The open sets $O_i$ in Theorem \ref{polypart} are called cells. Since they 
are disjoint and their union is $\mbb R^n \setminus Z(P)$, it follows that 
each cell is a union of connected components of $\mbb R^n \setminus Z(P)$. 
Moreover, since $Z(P)$ has Lebesgue measure zero, it follows from the 
inequalities in part (ii) of the theorem that
\begin{displaymath}
C_n^{-1} D^n \leq \# \{ O_i \} \leq C_n D^n.
\end{displaymath}
If $i \not= i'$, $x \in O_i$, and $y \in O_{i'}$, then $x$ and $y$ can not
belong to the same connected component of $\mbb R^n \setminus Z(P)$, and 
hence the line segment $\{ (1-t)x + ty : 0 \leq t \leq 1 \}$ must intersect
the zero set $Z(P)$. From this it follows that
\begin{itemize}
\item any line in $\mbb R^n$ can intersect at most $D+1$ of the cells $O_i$.
\end{itemize}

Let $2 < \alpha \leq 3$, $H$ be a weight in ${\mathcal F}_\alpha$, and $f$ a
function in $\Lambda(R,K,m)$. 

We apply Theorem \ref{polypart} in $\mbb R^3$ with the function $F$ defined 
by
\begin{displaymath} 
F(x) = \left\{ \begin{array}{ll}
         \mbox{Br}_\beta Ef(x)^p H(x) & \mbox{if $|x| \leq R$,} \\
		 0                            & \mbox{if $|x| > R$,}
		       \end{array} \right.
\end{displaymath}
and a degree $D$ that will be determined later in the argument.

In order for the resulting polynomial partitioning of $\mbb R^3$ to help us 
upgrade the two-dimensional bilinear restriction theorem into an estimate on
$\int F(x) dx$, we need the above property about intersections of lines and 
cells to also hold for the tubes $T \in \mbb T$ that are involved in the 
wave packet decomposition of $Ef_\tau$ as given by (\ref{finbuterr}) (with
$f$ replaced by $f_\tau$). But this is not true; in fact, it is possible for 
a tube $T$ to intersect all of the $\sim D^3$ cells $O_i$. 

Guth mended the situation by replacing the zero set $Z(P)$ by the {\it 
cell-wall} $W$ that he defined as the $R^{(1/2)+\delta}$-neighborhood of 
$Z(P)$, and the cells $O_i$ by the {\it modified cells} 
$O_i'= O_i \setminus W$. If one of our tubes $T$ intersects more than $D+1$ 
of the modified cells $O_i'$, then its core line intersects more than $D+1$ 
of the cells $O_i$, which implies that the core line lies in $Z(P)$. Since 
the radius of $T$ is $R^{(1/2)+\delta}$, it follows that $T$ lies in $W$, 
which is a contradiction. Therefore,
\begin{itemize}
\item any tube $T \in \mbb T$ can intersect at most $D+1$ of the modified 
      cells $O_i'$.
\end{itemize}

We write
\begin{displaymath} 
\int F(x) dx = \int_W F(x) dx + \int_{\cup_i O_i'} F(x) dx.
\end{displaymath}
If $\int_W F(x) dx \leq \int_{\cup_i O_i'} F(x) dx$, we say we are in the 
{\it cellular case}. Otherwise, we are in the {\it algebraic case}. 

In the cellular case, the above property about the intersection of tubes 
with modified cells will allow us to estimate our integral by induction on 
$\sum_\tau \| f_\tau \|_{L^2(\sigma)}^2$. The algebraic case is much more 
complicated.

In the algebraic case, the wave packets $Ef_{\tau,T}$ that contribute to our 
integral are those with `supporting' tubes $T$ that intersect the cell-wall 
$W$ in the ball $B(0,R)$. Guth separated such tubes into two groups 
depending on whether they lie in $B(0,R) \cap W$ for a significant portion 
of their length, or they cut $B(0,R) \cap W$ transversely. Here is the 
precise set-up.

We first recall that $W$ is the $R^{(1/2)+\delta}$-neighborhood of the zero 
set $Z(P)$ and we cover $B(0,R) \cap W$ by a finitely overlapping family 
$\{ B_j \}$ of balls each of radius $R^{1-\delta}$. Next, we let $Z_0(P)$ be 
the set of all non-singular points of $Z(P)$, and, for $z \in Z_0(P)$, we 
denote by $T_z Z(P)$ the tangent space to $Z(P)$ at the point $z$. Also, for 
each tube $T \in \mbb T$, we denote by $v(T)$ the unit vector in the 
direction of $T$.

Given a tube $T \in \mbb T$, we say $T$ cuts $W$ tangentially in $B_j$, and
write $T \in \mbb T_{j,\mbox{\rm \tiny tang}}$, if
\begin{displaymath}
\left\{ \begin{array}{l}
        T \cap W \cap B_j \not= \emptyset \\
        \mbox{Angle}(v(T),T_zZ(P)) \leq R^{-(1/2)+2\delta}
        \;\; \forall \; z \in Z_0(P) \cap 2B_j \cap 10T.
        \end{array} \right.
\end{displaymath}
Also, we say $T$ cuts $W$ transversely in $B_j$, and write 
$T \in \mbb T_{j,\mbox{\rm \tiny trans}}$, if
\begin{displaymath}
\left\{ \begin{array}{l}
        T \cap W \cap B_j \not= \emptyset \\
        \exists \; z \in Z_0(P) \cap 2B_j \cap 10T
        \mbox{ such that Angle}(v(T),T_zZ(P)) > R^{-(1/2)+2\delta}.
        \end{array} \right.
\end{displaymath}
We also let
\begin{displaymath}
f_{\tau,j,\mbox{\tiny tang}}
= \sum_{T \in \mbb T_{j,\mbox{\tiny tang}}} f_{\tau,T}
\hspace{0.25in} \mbox{and} \hspace{0.25in}
f_{j,\mbox{\tiny tang}}= \sum_\tau f_{\tau,j,\mbox{\tiny tang}},
\end{displaymath}
and
\begin{displaymath}
f_{\tau,j,\mbox{\tiny trans}}
= \sum_{T \in \mbb T_{j,\mbox{\tiny trans}}} f_{\tau,T}
\hspace{0.25in} \mbox{and} \hspace{0.25in}
f_{j,\mbox{\tiny trans}}= \sum_\tau f_{\tau,j,\mbox{\tiny trans}}.
\end{displaymath}

Recall that ${\mathcal P}$ is covered by $\sim K^2$ caps $\tau$ of diameter 
$\sim 1/K$. If $I$ is any subset of these caps, we write
\begin{displaymath}
f_{I,j,\mbox{\tiny trans}}= \sum_{\tau \in I} f_{\tau,j,\mbox{\tiny trans}}.
\end{displaymath}

We say that two caps $\tau_1$ and $\tau_2$ are non-adjacent if the distance
between them is $\geq 1/K$, and define
\begin{displaymath}
\mbox{Bil}_{P,\delta} E f_{j,\mbox{\tiny tang}} =
\sum_{\tau_1, \tau_2 \, \mbox{\tiny non-adjacent}}
|E f_{\tau_1,j,\mbox{\tiny tang}}|^{1/2}
|E f_{\tau_2,j,\mbox{\tiny tang}}|^{1/2}.
\end{displaymath}
We call the function $\mbox{Bil}_{P,\delta} E f_{j,\mbox{\tiny tang}}$ the
tangential part of $Ef$ with respect to the polynomial $P$ and the parameter
$\delta$.

The main properties of $\mbb T_{j,\mbox{\tiny tang}}$, 
$\mbb T_{j,\mbox{\tiny trans}}$, and 
$\mbox{Bil}_{P,\delta} E f_{j,\mbox{\tiny tang}}$ are presented in the 
following proposition.

\begin{alphprop}[\cite{guth:poly}]
\label{tangtrans}
{\rm (i)} For each $j$, 
$\mbb T_{j,\mbox{\tiny tang}} \cap \mbb T_{j,\mbox{\tiny trans}}=\emptyset$.
\\
{\rm (ii)} A tube $T \in \mbb T$ can belong to at most $D^{O(1)}$ different 
sets $\mbb T_{j,\mbox{\rm \tiny trans}}$.
\\
{\rm (iii)} For each $j$, the number of different $\theta$ such that
$\mbb T_{j,\mbox{\rm \tiny tang}} \cap \mbb T(\theta) \not= \emptyset$ is at
most $D^2 R^{(1/2)+O(\delta)}$.
\\
{\rm (iv)} Suppose $0 < \epsilon \leq 2$, $K \geq \sqrt[98]{10}$,
$K^{-\epsilon} \leq \beta \leq 1$, $\beta m \leq 10^{-5}$, and
$R \geq C K^\epsilon$. If $x \in B_j \cap W$ and $C$ is sufficiently large,
then
\begin{eqnarray*}
\mbox{\rm Br}_\beta E f(x)
& \leq & \frac{5}{4}
         \sum_I \mbox{\rm Br}_{2\beta} E f_{I,j,\mbox{\rm \tiny trans}}(x)
 + K^{100} \, \mbox{\rm Bil}_{P,\delta} E f_{j,\mbox{\rm \tiny tang}}(x) \\
& & + \; O \Big( R^{-N+1} \sum_\tau \| f_\tau \|_{L^1(\sigma)} \Big).
\end{eqnarray*}
\end{alphprop}

For a proof of part (i) of Proposition \ref{tangtrans}, we refer the reader 
to the paragraph immediately following 
\cite[Definitions 3.3 and 3.4]{guth:poly}. 

Part (ii) of Proposition \ref{tangtrans} is \cite[Lemma 3.5]{guth:poly}. 

Part (iii) of Proposition \ref{tangtrans} is \cite[Lemma 3.6]{guth:poly}. 
(For the very interesting connection between this result and the Kakeya 
problem, we also refer the reader to \cite[Conjecture B.1]{guth:poly2} and 
\cite{katzrogers}.)

Part (iv) is \cite[Lemma 3.8]{guth:poly}. (For a proof of this part in the 
precise way it is stated in Proposition \ref{tangtrans}, we refer to   
\cite[Lemma 9-B]{plms12046}).)

Part (iv) allows us to bound $\int_{B(0,R)} F(x) dx$ in the algebraic case 
by induction over the radius $R$ (recall that each $B_j$ has radius 
$R^{1-\delta}$) provided we have a bound on
\begin{displaymath}
\int_{B_j \cap W} \mbox{\rm Bil}_{P,\delta} 
E f_{j,\mbox{\rm \tiny tang}}(x)^p H(x) dx.
\end{displaymath}

The above discussion was formulated in \cite{plms12046}) in the following 
theorem (which is in turn a reformulation of \cite[Theorem 3.1]{guth:poly} 
in the weighted setting). We remind the reader that the space 
$\Lambda(R,K,m)$ of functions was defined in the paragraph before the last 
of Subsection 4.2.  

\begin{alphthm}[{\cite[Theorem 9.1]{plms12046})}]
\label{biltobr}
Let $3< p \leq 4$, $\epsilon > 0$, $0 \leq q_1 \leq 1 \leq 2 q_0$, 
$0< q_2 < q_0$, and $H$ be a weight of fractal dimension $\alpha$. Also, let 
$\delta=\epsilon^2$, $\delta_{\mbox{\rm \tiny deg}}= \epsilon^4$, and 
$\delta_{\mbox{\rm \tiny trans}}= \epsilon^6$.

Suppose that
\begin{eqnarray}
\label{estonbil}
\lefteqn{\int_{B_j \cap W}
\mbox{\rm Bil}_{P,\delta} E f_{j,\mbox{\rm \tiny tang}}(x)^p H(x) dx}
\nonumber \\
& \leq & C_{\epsilon,K} R^{O(\delta)} R^{q_2\epsilon} A_\alpha(H)^{q_1}
        \Big( \sum_\tau \| f_\tau \|_{L^2(S)}^2 \Big)^{(3/2)+\epsilon}
\end{eqnarray}
whenever $R \geq C$, $K \geq 100$, $m \geq 1$, $f \in \Lambda(R,K,m)$, and 
$P \in \mbox{\rm Poly}_D(\mbb R^3)$ with 
$D=R^{\delta_{\mbox{\rm \tiny deg}}}$ and $P$ a product of non-singular 
polynomials.

Then there is a constant $c_0$, which is independent of $q_0$, $q_1$, and 
$p$, such that if $\epsilon \leq \min[ c_0,(p-3)/2]$, then there is a
$K=K(\epsilon)$ such that
\begin{eqnarray}
\label{underbil}
\lefteqn{\int_{B(0,R)} \mbox{\rm Br}_\beta E f(x)^p H(x) dx} \nonumber \\
& \leq & C_\epsilon R^{q_0\epsilon} A_\alpha(H)^{q_1}
         \Big( \sum_\tau \| f_\tau \|_{L^2(S)}^2 \Big)^{(3/2)+\epsilon}
         R^{\delta_{\mbox{\rm \tiny trans}} \log(K^\epsilon \beta m)}
\end{eqnarray}
for all $\beta \geq K^{-\epsilon}$, $m \geq 1$, $R \geq 1$, and
$f \in \Lambda(R,K,m)$. Moreover,
$\lim_{\epsilon \to 0} K(\epsilon)=~\!\!\infty$.
\end{alphthm}

In \cite[Theorem 9.1]{plms12046}), there is an additional parameter $b$. 
Theorem \ref{biltobr} is the special case $b=1$.

\section{Proof of Theorem \ref{mainjj}}

There are two known estimates that are key to proving Theorem \ref{mainjj}.
The first, which is part (i) of the next lemma, was proved in 
\cite{guth:poly} by adapting Cordoba's $L^4$ argument from \cite{lfour} to
neighborhoods of algebraic surfaces in $\mbb R^3$. The second, which is part 
(ii) of the next lemma, was proved in \cite{dgowwz:falconer} by adapting the 
refined bilinear Strichartz estimates from \cite{dgl:schrodinger} to the 
weighted restriction setting.

\begin{alphlemma}
\label{2knownest}
Suppose $B_j$ is one of the $R^{1-\delta}$-balls that cover $B(0,R) \cap W$, 
and $\tau_1$ and $\tau_2$ are non-adjacent $K^{-1}$-caps. Then:
\\
{\rm (i)} We have
\begin{eqnarray*}
\lefteqn{\int_{B_j \cap W} |E f_{\tau_1,j,\mbox{\tiny tang}}|^2
         |E f_{\tau_2,j,\mbox{\tiny tang}}|^2 dx} \nonumber \\
& \lct & R^{O(\delta)} R^{-1/2} 
         \| f_{\tau_1,j,\mbox{\tiny tang}} \|_{L^2(\sigma)}^2
		 \| f_{\tau_2,j,\mbox{\tiny tang}} \|_{L^2(\sigma)}^2
         + \mbox{\rm negligible}
\end{eqnarray*}
for all $f \in \Lambda(R,K,m)$.
\\
{\rm (ii)} We have
\begin{displaymath}
\int_{B_j} |E f_{\tau_1,j,\mbox{\tiny tang}}|^2
           |E f_{\tau_2,j,\mbox{\tiny tang}}|^2 H dx 
\lct R^{O(\delta)} A_\alpha(H) R^{(\alpha-2)/4} 
     \| f_{\tau_1} \|_{L^2(\sigma)}^{3/2} 
	 \| f_{\tau_2} \|_{L^2(\sigma)}^{3/2} 
\end{displaymath}
for all $f \in \Lambda(R,K,m)$.
\end{alphlemma}

Part (i) of Lemma \ref{2knownest} is the estimate immediately preceding 
\cite[Equation (43)]{guth:poly}. 

Part (ii) of Lemma \ref{2knownest} is 
\cite[Equation (4.19)]{dgowwz:falconer}. We alert the reader, however, that
\cite[Equation (4.19)]{dgowwz:falconer} is stated as
\begin{displaymath}
\int_{B_j} |E f_{\tau_1,j,\mbox{\tiny tang}}|^2
           |E f_{\tau_2,j,\mbox{\tiny tang}}|^2 {\mathcal H} dx 
\lct R^{O(\delta)} R^{3 \gamma_3^0} \| f_{\tau_1} \|_{L^2(\sigma)}^{3/2} 
                                    \| f_{\tau_2} \|_{L^2(\sigma)}^{3/2} 
\end{displaymath} 
with $3 \gamma_3^0 = \alpha-2$ and 
${\mathcal H} \in {\mathcal F}_{\alpha,3}$, but its proof gives the better 
exponent $(\alpha-2)/4$ of $R$ that appears in the above lemma (see the last 
line of \cite[Section 4]{dgowwz:falconer}). Also, the equivalence between 
(\ref{unify1}) and (\ref{unify2}) allows us to replace the
${\mathcal H} \in {\mathcal F}_{\alpha,3}$ by a weight $H$ of fractal 
dimension $\alpha$.

To prove our theorem, we use Lemma \ref{2knownest} to establish the estimate
(\ref{estonbil}) on the tangential part of $Ef$, which we then turn into an 
estimate on the broad part of $Ef$ via Theorem \ref{biltobr}. We then 
localize the weight function as described in Section 1 and use parabolic 
scaling and induction on the radius to turn the estimate on the broad part 
of $Ef$ into the desired estimate of Theorem \ref{mainjj}.

\subsection{Estimate on the tangential part}

The function we are interested in estimating is
\begin{displaymath}
\mbox{Bil}_{P,\delta} E f_{j,\mbox{\tiny tang}} =
\sum_{\tau_1, \tau_2 \, \mbox{\tiny non-adjacent}}
|E f_{\tau_1,j,\mbox{\tiny tang}}|^{1/2}
|E f_{\tau_2,j,\mbox{\tiny tang}}|^{1/2}.
\end{displaymath}
To simplify the notation a little, we write
\begin{displaymath}
F = \mbox{Bil}_{P,\delta} E f_{j,\mbox{\tiny tang}} 
\hspace{0.25in} \mbox{ and } \hspace{0.25in}
G= |E f_{\tau_1,j,\mbox{\tiny tang}}| \, |E f_{\tau_2,j,\mbox{\tiny tang}}|.
\end{displaymath}
Since the constant $C_{\epsilon,K}$ in (\ref{estonbil}) in Theorem 
\ref{biltobr} is allowed to depend on $K$, and there are $\sim K^4$ pairs 
$\tau_1, \tau_2$, we have
\begin{displaymath}
\int_{B_j \cap W} F^p H dx 
\lct \sum_{\tau_1, \tau_2 \, \mbox{\tiny non-adjacent}}
\int_{B_j \cap W} G^{p/2} H dx.
\end{displaymath}

Our plan is to apply H\"{o}lder's inequality to 
$\int_{B_j \cap W} G^{p/2} H dx$ and then use both parts of Lemma 
\ref{2knownest}. In order to do this, we need to find positive numbers $a$ 
and $b$ such that $a+b=p/2$, $aq=2$, and $bq'=3/2$ for some $1< q < \infty$. 
This means $a$ and $b$ must satisfy
\begin{displaymath}
\left\{ \begin{array}{lll}
        \frac{a}{2}+\frac{2b}{3} & = & 1 \\
		a+b                      & = & \frac{p}{2}
		\end{array} \right.
\end{displaymath}
which gives $a=2(p-3)$ and $b=(3/2)(4-p)$, and hence
\begin{displaymath}
\int_{B_j \cap W} G^{p/2} H dx =
\int_{B_j \cap W} G^{2(p-3)} G^{(3/2)(4-p)} H dx.
\end{displaymath}
Applying H\"{o}lder's inequality with $q=2/a=1/(p-3)$, we get
\begin{displaymath}
\int_{B_j \cap W} F^p H dx 
\lct \sum_{\tau_1, \tau_2 \, \mbox{\tiny non-adjacent}}
\Big( \int_{B_j \cap W} G^2 H dx \Big)^{p-3}
\Big( \int_{B_j \cap W} G^{3/2} H dx \Big)^{4-p}.
\end{displaymath}
Writing
\begin{displaymath}
G_0= \| f_{\tau_1,j,\mbox{\tiny tang}} \|_{L^2(\sigma)} 
     \| f_{\tau_2,j,\mbox{\tiny tang}} \|_{L^2(\sigma)},
\end{displaymath}
Lemma \ref{2knownest} tells us that
\begin{displaymath}
\int_{B_j \cap W} G^2 H dx \lct R^{O(\delta)} R^{-1/2} G_0^2
\end{displaymath}
(because $\| H \|_{L^\infty(\mbb R^3)} \leq 1$) and
\begin{displaymath}
\int_{B_j \cap W} G^{3/2} H dx 
\lct R^{O(\delta)} A_\alpha(H) R^{(\alpha-2)/4} G_0^{3/2}.
\end{displaymath}
Therefore,
\begin{equation}
\label{bding0}
\int_{B_j \cap W} F^p H dx \lct R^{O(\delta)} A_\alpha(H)^{4-p} R^{-(p-3)/2}
R^{(\alpha-2)(4-p)/4} \sum_{\tau_1, \tau_2} G_0^{p/2}.
\end{equation}

Following \cite{guth:poly}, we now use the fact that only few of the caps
$\theta$ contribute to $G_0$. 

Suppose $\tau$ is one of our $K^{-1}$-caps, and let 
\begin{displaymath}
\{ \theta_1, \ldots, \theta_L \} 
= \{ \theta : \theta \cap \tau \not= \emptyset \mbox{ and }
\mbb T_{j,\mbox{\tiny tang}} \cap \mbb T(\theta) \not= \emptyset \}
\end{displaymath}
and ${\mathcal T}_l = \mbb T_{j,\mbox{\tiny tang}} \cap \mbb T(\theta_l)$.
Then
\begin{displaymath}
f_{\tau,j,\mbox{\tiny tang}} 
= \sum_{l=1}^L \sum_{T \in {\mathcal T}_l} f_{\tau,T}
\end{displaymath}
and it follows by (\ref{fewcaps}) that
\begin{displaymath}
\int_{3\theta_l} |f_{\tau,j,\mbox{\tiny tang}}|^2 d\sigma 
\lct \int_{10\theta_l} |f_\tau|^2 d\sigma
\end{displaymath} 
for all $1 \leq l \leq L$. By part (i) of Proposition \ref{wave}, we know 
that $f_{\tau,j,\mbox{\tiny tang}}$ is supported in 
$\cup_{l=1}^L 3 \theta_l$, so
\begin{displaymath}
\int |f_{\tau,j,\mbox{\tiny tang}}|^2 d\sigma 
\lct \sum_{l=1}^L \int_{10\theta_l} |f_\tau|^2 d\sigma.
\end{displaymath}
On the one hand, this gives
\begin{equation}
\label{onehand}
\int |f_{\tau,j,\mbox{\tiny tang}}|^2 d\sigma \lct \int |f_\tau|^2 d\sigma.
\end{equation}
On the other hand, from the definition of $\Lambda(R,K,m)$ (see
(\ref{baverages})) we know that
\begin{displaymath}
\int_{10\theta_l} |f_\tau|^2 d\sigma \lct \frac{1}{R},
\end{displaymath} 
so
\begin{displaymath}
\int |f_{\tau,j,\mbox{\tiny tang}}|^2 d\sigma 
\lct \sum_{l=1}^L \int_{10\theta_l} |f_\tau|^2 d\sigma \lct \frac{L}{R}.
\end{displaymath}
Part (iii) of Proposition \ref{tangtrans} also tells us that 
$L \lct R^{O(\delta)} R^{1/2}$, so
\begin{displaymath}
\int |f_{\tau,j,\mbox{\tiny tang}}|^2 d\sigma \lct R^{O(\delta)} R^{-1/2}.
\end{displaymath}
Therefore,
\begin{displaymath}
G_0 = \| f_{\tau_1,j,\mbox{\tiny tang}} \|_{L^2(\sigma)} 
      \| f_{\tau_2,j,\mbox{\tiny tang}} \|_{L^2(\sigma)} 
\lct R^{O(\delta)} R^{-1/2}.
\end{displaymath}

Going back to (\ref{bding0}), we now have
\begin{displaymath}
G_0^{p/2}= G_0^{(p-3-2\epsilon)/2} G_0^{(3/2)+\epsilon}
\lct R^{O(\delta)} R^{-(p-3-2\epsilon)/4} G_0^{(3/2)+\epsilon},
\end{displaymath}
and (\ref{bding0}) becomes
\begin{displaymath}
\int_{B_j \cap W} F^p H dx 
\lct R^{O(\delta)} R^{\epsilon/2} A_\alpha(H)^{4-p} R^{-3(p-3)/4}
R^{(\alpha-2)(4-p)/4} \sum_{\tau_1, \tau_2} G_0^{(3/2)+\epsilon}.
\end{displaymath}
By easy considerations,
\begin{displaymath}
\sum_{\tau_1, \tau_2} G_0^{(3/2)+\epsilon}= \Big( \sum_\tau 
\| f_{\tau,j,\mbox{\tiny tang}} \|_{L^2(\sigma)}^{(3/2)+\epsilon} \Big)^2
\lct \Big( \sum_\tau 
\| f_{\tau,j,\mbox{\tiny tang}} \|_{L^2(\sigma)}^2 \Big)^{(3/2)+\epsilon},
\end{displaymath}
and so, by (\ref{onehand}),
\begin{displaymath}
\sum_{\tau_1, \tau_2} G_0^{(3/2)+\epsilon}
\lct \Big( \sum_\tau \| f_\tau \|_{L^2(\sigma)}^2 \Big)^{(3/2)+\epsilon}.
\end{displaymath}
Also,
\begin{displaymath}
R^{-3(p-3)/4} R^{(\alpha-2)(4-p)/4} = 1 
\hspace{0.25in} \mbox{ if } \hspace{0.25in} p = \frac{4\alpha+1}{\alpha+1}.
\end{displaymath}
Therefore, with this value of $p$, we have
\begin{displaymath}
\int_{B_j \cap W} F^p H dx 
\lct R^{O(\delta)} R^{\epsilon/2} A_\alpha(H)^{4-p}
\Big( \sum_\tau \| f_\tau \|_{L^2(\sigma)}^2 \Big)^{(3/2)+\epsilon}.
\end{displaymath}

\subsection{Estimate on the broad part}

In the previous subsection, we established (\ref{estonbil}) with 
$p=(4\alpha+1)/(\alpha+1)$, $q_1=4-p$, and $q_2=1/2$. Inserting these values
into Theorem \ref{biltobr} and choosing $q_0=1$, we conclude that there is a
constant $c_0$ such that to every 
$\epsilon \leq \min[c_0,(\alpha-2)/(2\alpha+2)]$ there are constants 
$K=K(\epsilon)$ and $C_\epsilon$ such that 
$\lim_{\epsilon \to 0} K(\epsilon)= \infty$ and
\begin{eqnarray*}
\lefteqn{\int_{B(0,R)} \mbox{Br}_{\beta} E f(x)^p H(x) dx} \\
& \leq & C_\epsilon \, R^\epsilon A_\alpha(H)^{3/(\alpha+1)}
         \Big( \sum_\tau \| f_\tau \|_{L^2(\sigma)}^2 \Big)^{(3/2)+\epsilon}
         R^{\delta_{\mbox{\tiny trans}} \log(K^\epsilon \beta m)}
\end{eqnarray*}
for all $\beta \geq K^{-\epsilon}$, $m \geq 1$, $R \geq 1$, and
$f \in \Lambda(R,K,m)$.

Suppose $f \in L^\infty(\sigma)$ is a non-zero function, and define the 
function $g : {\mathcal P} \to \mbb C$ by 
$g=\pi^{-1/2} \| f \|_{L^\infty(\sigma)}^{-1} f$, so that 
\begin{displaymath}
\int_{B(\xi_0,R^{-1/2}) \cap {\mathcal P}} |g|^2 d\sigma \leq \frac{1}{R}
\end{displaymath}
for all $\xi_0 \in {\mathcal P}$. Writing $g=\sum_\tau g_\tau$ with
$\mbox{supp} \, g_\tau \subset \tau$ and
$(\mbox{supp} \, g_\tau) \cap (\mbox{supp} \, g_{\tau'})= \emptyset$ if
$\tau \not= \tau'$, we see that $g \in \Lambda(R,K,m)$ and
$\sum_\tau \| g_\tau \|_{L^2(\sigma)}^2= \| g \|_{L^2(\sigma)}^2$. Applying 
the above estimate on the broad part with $\beta= K^{-\epsilon}$ and $m=1$, 
we obtain
\begin{displaymath}
\int_{B(0,R)} \mbox{Br}_{K^{-\epsilon}} E g(x)^p H(x) dx
\leq C_\epsilon \, R^\epsilon A_\alpha(H)^{3/(\alpha+1)} 
     \| g \|_{L^2(\sigma)}^{3+2\epsilon}.
\end{displaymath}
Since $\| g \|_{L^2(S)} \leq 1$, we have 
$\| g \|_{L^2(\sigma)}^{3+2\epsilon} \leq \| g \|_{L^2(\sigma)}^3$. So,
replacing $g$ by $f$, and letting $c=\min[c_0,(\alpha-2)/(2\alpha+2)]$, we 
see that to every $\epsilon \leq c$ there are constants 
$K=K(\epsilon)$ and $\bar{C}_\epsilon$ such that 
$\lim_{\epsilon \to 0} K(\epsilon)= \infty$ and
\begin{equation}
\label{priortolocal}
\int_{B(0,R)} \mbox{Br}_{K^{-\epsilon}} E f(x)^p H(x) dx
\leq \tilde{C}_\epsilon \, R^\epsilon A_\alpha(H)^{3/(\alpha+1)}
     \| f \|_{L^2(\sigma)}^3 \| f \|_{L^\infty(\sigma)}^{p-3}    
\end{equation}
for all functions $f \in L^\infty(\sigma)$, weights $H$ of fractal dimension
$\alpha$, and radii $R \geq 1$.

\subsection{Estimate on $Ef$}

Moving from the estimate on the $L^p$-norm of 
$\mbox{Br}_{K^{-\epsilon}} E f(x)$ against the weight $H$ into an estimate 
on the $L^p$-norm of $Ef$ itself involves localizing the weight $H$ at scale
$r = K^{-1}$, followed by parabolic scaling, followed by an induction on the 
radius argument. These three arguments were carried out in detail by the 
author in \cite{plms12046} and resulted in the following theorem that we now
borrow from that paper.

\begin{alphthm}[{\cite[Theorem 12.1]{plms12046})}]
\label{parabscaling}
Let $0 < \alpha \leq 3$, $3 \leq p \leq 4$, $2 \leq \gamma \leq 3$,
$0 \leq q_1 \leq 1$, and $c > 0$.

Suppose that we have the following estimate on the broad part of $Ef$: to
every $\epsilon \in (0,c)$ there are constants $K(\epsilon)$ and
$\bar{C}_\epsilon$ such that $\lim_{\epsilon \to 0} K(\epsilon)= \infty$ and
\begin{equation}
\label{parabscalcond}
\int_{B(0,R)} \mbox{\rm Br}_{K^{-\epsilon}} Ef(x)^p H(x)dx
\leq \bar{C}_\epsilon R^\epsilon A_\alpha(H)^{q_1}
     \| f \|_{L^2(\sigma)}^{\gamma} \| f \|_{L^\infty(\sigma)}^{p-\gamma}
\end{equation}
for all radii $R \geq 1$, weights $H$ of dimension $\alpha$, and functions 
$f \in L^\infty(\sigma)$.

If $2p-\alpha-1-\gamma > 0$, then there is a constant $c'$, which only
depends on $\alpha, p$, and $\gamma$, such that for $0 < \epsilon < c'$ we
have
\begin{displaymath}
\int_{B(0,R)} |Ef(x)|^p H(x)dx \leq C_\epsilon R^\epsilon
\Big( \max \big[ A_\alpha(H), A_\alpha(H)^{q_1} \big] \Big)
\| f \|_{L^2(\sigma)}^{\gamma} \| f \|_{L^\infty(\sigma)}^{p-\gamma},
\end{displaymath}
with
\begin{displaymath}
C_\epsilon=
2 \big( \bar{C}_\epsilon + 10^4 \sigma(S)^4 \big),
\end{displaymath}
for all radii $R \geq 1$, weights $H$ of dimension $\alpha$, and functions 
$f \in L^\infty(\sigma)$.
\end{alphthm}

Since
\begin{displaymath}
\| f \|_{L^2(\sigma)}^3
= \| f \|_{L^2(\sigma)}^{\gamma} \| f \|_{L^2(\sigma)}^{3-\gamma}
\lct \| f \|_{L^2(\sigma)}^{\gamma} \| f \|_{L^\infty(\sigma)}^{3-\gamma},
\end{displaymath}
the estimate (\ref{priortolocal}) implies (\ref{parabscalcond}) with 
$p=(4\alpha+1)/(\alpha+1)$ and $q_1=4-p$, and hence Theorem 
\ref{parabscaling} tells us that the estimate
\begin{displaymath}
\int_{B(0,R)} |Ef(x)|^p H(x)dx \leq C_\epsilon R^\epsilon
\Big( \max \big[ A_\alpha(H), A_\alpha(H)^{q_1} \big] \Big)
\| f \|_{L^2(\sigma)}^{\gamma} \| f \|_{L^\infty(\sigma)}^{p-\gamma}
\end{displaymath}
holds whenever $2p-\alpha-1-\gamma > 0$. So, to finish the proof of Theorem
\ref{mainjj}, we just have to show that this estimate remains true if one 
replaces $\max[A_\alpha(H),A_\alpha(H)^{q_1}]$ by $A_\alpha(H)$:
\begin{displaymath}
\int_{B(0,R)} |Ef(x)|^p H(x) dx 
\leq C_\epsilon R^\epsilon A_\alpha(H) 
\| f \|_{L^2(\sigma)}^\gamma \| f \|_{L^\infty(\sigma)}^{p-\gamma}. 
\end{displaymath}
In fact, since the unit paraboloid ${\mathcal P}$ is compact, we can find a 
$C_0^\infty$ function $\phi$ on $\mbb R^3$ that satisfies $|\phi| \geq 1$ on 
${\mathcal P}$ and $\widehat{\phi}$ is supported in the unit ball. So, if we 
define the function $g$ on ${\mathcal P}$ by $g=f/\phi$ and observe that 
$|g| \leq |f|$, $E f = (E g) \ast \widehat{\phi}$, and 
$|E f|^p \lct |E g|^p \ast |\widehat{\phi}|$, we see that 
\begin{eqnarray*}
\int |E f(x)|^p H(x) dx
& \lct & \int |E g(y)|^p \int |\widehat{\phi}(x-y)| H(x) dx dy \\ 
& \lct & A_\alpha(H) \int |E g(y)|^p \widetilde{H}(y) dy \\
& \lct & A_\alpha(H) R^\epsilon
         \max[A_\alpha(\widetilde{H}), A_\alpha(\widetilde{H})^{q_0}] 
	\| g \|_{L^2(\sigma)}^\gamma \| g \|_{L^\infty(\sigma)}^{p-\gamma} \\
& \lct & A_\alpha(H) R^\epsilon
         \max[A_\alpha(\widetilde{H}), A_\alpha(\widetilde{H})^{q_0}] 
	     \| f \|_{L^2(\sigma)}^\gamma \| f \|_{L^\infty(\sigma)}^{p-\gamma},
\end{eqnarray*}
where 
\begin{displaymath}
\widetilde{H}(y) 
= A_\alpha(H)^{-1} \| \widehat{\phi} \|_{L^\infty(\mbb R^3)}^{-1} 
  \int |\widehat{\phi}(x-y)| H(x) dx 
\end{displaymath}
is clearly a weight of fractal dimension $\alpha$ with 
$A_\alpha(\widetilde{H}) \lct 1$. Therefore,
\begin{displaymath}
\int_{B(0,R)} |Ef(x)|^p H(x) dx 
\leq C_\epsilon R^\epsilon A_\alpha(H) 
\| f \|_{L^2(\sigma)}^\gamma \| f \|_{L^\infty(\sigma)}^{p-\gamma}
\end{displaymath}
whenever $p=(4\alpha+1)/(\alpha+1)$ and $2p-\alpha-1-\gamma > 0$, as 
desired.

\end{document}